\documentclass[conference]{IEEEtran}
\IEEEoverridecommandlockouts

\usepackage{cite}
\usepackage{amsmath,amssymb,amsfonts}
\usepackage{algorithmic}
\usepackage{graphicx}
\usepackage{textcomp}
\usepackage{xcolor}
\def\BibTeX{{\rm B\kern-.05em{\sc i\kern-.025em b}\kern-.08em
    T\kern-.1667em\lower.7ex\hbox{E}\kern-.125emX}}

\newtheorem{teo}{Theorem}[section]

\newtheorem{prop}[teo]{Proposition}

\newtheorem{remark}{Remark}

\newcommand{\R}{\mathbb{R}} 
\newcommand{\W}{\mathcal{W}} 
\newcommand{\C}{\mathbb{C}} 

\begin{document}

\title{An optimal estimate for the norm of wavelet localization operators}

\author{\IEEEauthorblockN{Federico Riccardi}
\IEEEauthorblockA{\textit{Dipartimento di Scienze Matematiche ``Giuseppe Luigi Lagrange''} \\
\textit{Politecnico di Torino}\\
Torino, Italy \\
federico.riccardi@polito.it}
}

\maketitle

\begin{abstract}
In this paper we prove an optimal estimate for the norm of wavelet localization operators with Cauchy wavelet and weight functions that satisfy two constraints on different Lebesgue norms. We prove that multiple regimes arise according to the ratio of these norms: if this ratio belongs to a fixed interval (which depends on the Lebesgue exponents) then both constraints are active, while outside this interval one of the constraint is inactive. Furthermore, we characterize optimal weight functions.
\end{abstract}


\section{Introduction}
In \cite{nicolatilli_norm}, sharp estimates for the norm of time-frequency and wavelet localization operators were obtained using the remarkable Faber-Krahn inequality for the short-time Fourier transform with Gaussian window, proved in \cite{nicolatilli_fk}, and the analogous result for the wavelet transform with Cauchy wavelets proved in \cite{ramos-tilli}. More precisely, the problem of maximizing the operator norm of a (time-frequency or wavelet) localization operator under the assumption that the weight weight function belongs to some Lebesgue space for $1 \leq p < \infty$ (and possibly also in $L^{\infty})$ was addressed, along with the characterization of optimal weights. This result, for time-frequency localization operator, was generalized by the author in \cite{riccardi_new_optimal}, where weight functions are assumed to belong to the intersection of two Lebesgue spaces $L^p \cap L^q$, for $1<p,q<\infty$. In this paper, we address the analogous problem for wavelet localization operator, thus proving that similar phenomena appear in this setting.

\section{Preliminaries, notation and main result}

We consider the Fourier transform with the following normalization
\begin{equation*}
	\hat{f}(\omega) = \dfrac{1}{\sqrt{2\pi}} \int_{\R} f(t) e^{-i \omega t} \, dt.
\end{equation*}
For $\beta > 0$, we consider the so-called Cauchy wavelet $\psi_{\beta} \in L^2(\R)$ defined by
\begin{equation*}
	\widehat{\psi_{\beta}}(\omega) = c_{\beta} \chi_{[0,+\infty)}(\omega) \omega^{\beta} e^{-\omega}, \quad \omega \in \R,
\end{equation*}
where $c_{\beta}$ is chosen so that $ 2 \pi \|\widehat{\psi_{\beta}}\|_{L^2(\R_+, d\omega/\omega)} = 1$ ($\R_+ = (0,+\infty)$). The wavelet transform with respect to a Cauchy wavelet is then defined as
\begin{equation*}\label{eq : def wavelet transform}
	\W_{\psi_{\beta}} f (x,y) = \dfrac{1}{\sqrt{y}} \int_{\R} f(t) \overline{\psi_{\beta}\left(\dfrac{t-x}{y}\right)} \, dt,
\end{equation*}
where $(x,y) \in \R \times \R_+$. The wavelet transform is well defined from the Hardy space $H^2(\R) = \{f \in L^2(\R) \colon \allowbreak \hat{f}(\omega) \text{ for a.e. } \omega < 0\}$, with the norm inherited from $L^2(\R)$, into $L^2(\R \times \R_+, d\nu)$, where $d\nu = dx dy /y^2 $ is the left Haar measure on $\R \times \R_+ \simeq \C_+$ regarded as the ``$ax + b$'' group. The choice of the constant $c_{\beta}$ is such that $\W_{\psi_{\beta}} \colon H^2(\R) \rightarrow L^2(\R \times \R_+, d\nu)$ is an isometry.

Given $F \in L^p(\R \times \R_+, d\nu)$, $1 \leq p \leq \infty$, and $f,g \in H^2(\R)$, the wavelet localization operator $L_{F,\beta}$ with weight $F$ is defined through the duality
\begin{equation*}
	\langle L_{F, \beta} f, g \rangle_{L^2} = \int_{\R \times \R^+} F \W_{\psi_{\beta}} f \overline{\W_{\psi_{\beta}} g} \, d\nu
\end{equation*}
and is a bounded operator from $H^2(\R)$ into itself (see \cite{wong}). Our aim is to find optimal estimates
\begin{equation*}
	\|L_{F,\beta}\|_{H^2 \rightarrow H^2} \leq C(p,q,A,B)
\end{equation*}
where $F \in L^p(\C_+, d\nu) \cap L^q(\C_+,d\nu)$ satisfies the constraints
\begin{equation}\label{eq : constraints wavelet case}
	\|F\|_{L^p(\C_+, d\nu)} \leq A, \quad \|F\|_{L^q(\C_+, d\nu)} \leq B,
\end{equation}
for some $p,q \in (1,+\infty)$ (with $p \neq q$) and some $A,B \in (0,+\infty)$. For the ease of notation, from now on we denote the operator norm of $L_{F,\beta}$ simply as $\|L_{F,\beta}\|$.

Our proof, relies on the following theorem. 
\begin{teo}[\cite{nicolatilli_norm} Theorem 5.3]\label{th:estimate norm}
	Let $F \in L^p(\C_+, d\nu)$ for some $p \geq 1$ and let $v(t) = \nu(\{|F|>t\})$ be the distribution function of $|F|$. Then
	\begin{equation}\label{eq:estimate norm}
		\|L_{F,\beta}\| \leq \int_0^{\infty} G(v(t))\, dt,
	\end{equation}
	where $G(s) = 1 - (1+s/4\pi)^{-2\beta}$. Equality is achieved if and only if $F(z) = e^{i \theta} \rho(|z-z_0|^2/|z-\overline{z_0}|^2)$, for some $\theta \in \R$, $z_0 \in \C_+$ and some nonincreasing function $\rho \colon [0,1) \to [0,+\infty)$.
\end{teo}

Before stating our main theorem, we introduce the following constants
\begin{equation*}
	\alpha_p = \dfrac{p-1}{2\beta+1}, \quad \sigma_p = \dfrac{p-1}{2\beta p+1}, \quad \kappa_p = \dfrac{p-1}{p},
\end{equation*}
and the same for $q$.
\begin{teo}\label{th:main theorem wavelet case}
	Let $p, q \in (1,+\infty)$, $A,B \in (0,+\infty)$ and
	\begin{align*}
		r_1 &= (4 \pi)^{\frac{1}{q} - \frac{1}{p}} \left(\dfrac{\alpha_q}{p-\alpha_q}\right)^{-\frac{1}{p}} (\sigma_q)^{\frac{1}{q}},\\
		r_2 &= (4 \pi)^{\frac{1}{q} - \frac{1}{p}} \left(\dfrac{\alpha_p}{q-\alpha_p}\right)^{\frac{1}{q}} (\sigma_p)^{-\frac{1}{p}}.
	\end{align*}
	Let $ F \in L^p(\C_+, d\nu) \cap L^q(\C_+, d\nu)$ satisfying \eqref{eq : constraints wavelet case}. Then:
	\begin{enumerate}
		\item If $q > \alpha_p$ and $B/A \geq r_2$ (respectively $p > \alpha_q$ and $B/A \leq r_1$) then the $L^q$ (respectively $L^p$) constraint is inactive and it holds
		\begin{align*}\label{eq:estimate wavelet first case}
			\|L_{F,\beta} \| &\leq \dfrac{2\beta}{(4\pi)^{1/p}} \sigma_p^{\kappa_p} A\\
			(\text{respectively } \|L_{F,\beta} \| &\leq \dfrac{2\beta}{(4\pi)^{1/q}} \sigma_q^{\kappa_q} B )
		\end{align*}
		with equality if and only if, for some $\theta \in \R$ and some $z_0 \in \C_+$, it holds
		\begin{align*}
			F(z) &= e^{i \theta} \lambda \left(1 - \left| \tfrac{z-z_0}{z-\overline{z_0}} \right|^2 \right)^{-\alpha_p} \\ (\text{respectively } F(z) &= e^{i \theta} \lambda \left(1 - \left| \tfrac{z-z_0}{z-\overline{z_0}} \right|^2 \right)^{-\alpha_q}),
		\end{align*}
		with $z \in \C_+$ and $\lambda = A(4 \pi \sigma_p)^{-1/p}$ (respectively $\lambda = B(4 \pi \sigma_q)^{-1/q}$).
		\item Otherwise it holds
		\begin{equation}\label{eq:estimate wavelet second case}
			\|L_{F,\beta} \| \leq \int_{0}^{\infty} G(u(t)) \, dt,
		\end{equation}
		where
		\begin{equation}\label{eq:expression of maximizer distribution function}
			u(t) = 4 \pi \max \{ (\lambda_1 t^{p-1} + \lambda_2 t^{q-1})^{-\frac{1}{2\beta+1}}-1, \, 0\}
		\end{equation}
		and $\lambda_1, \lambda_2>0$ are uniquely determined by
		\begin{equation*}\label{eq:equations for multipliers hyperbolic}
			\begin{cases}
				p \int_0^{+\infty} t^{p-1} u(t) \,dt = A^p\\
				q \int_0^{+\infty} t^{q-1} u(t) \,dt = B^q
			\end{cases}.
		\end{equation*}
		Moreover, letting $T > 0$ denote the unique solution of $\lambda_1 T^{p-1} + \lambda_2 T^{q-1} = 1$, the function $t \mapsto \{(\lambda_1 t^{p-1} + \lambda_2 t^{q-1})^{-\frac{1}{2\beta+1}}-1\}$ defined on $(0,T]$ is invertible and we denote by $\psi : [0, + \infty) \rightarrow (0,T]$ its inverse. Then, equality in \eqref{eq:estimate wavelet second case} is achieved if and only if
		\begin{equation*}
			F(z) = e^{i\theta} \psi\left(\frac{|z-z_0|^2/|z-\overline{z_0}|^2}{1-|z-z_0|^2/|z-\overline{z_0}|^2}\right)
		\end{equation*}
		for some $\theta \in \R$ and some $z_0 \in \C_+$.
	\end{enumerate}
\end{teo}

\begin{IEEEproof}[Proof of \ref{th:main theorem wavelet case}(1)]
	The proof of the first part of the theorem essentially follows from the observation that optimal functions for the problem with just the $L^p$ constraint (see \cite[Theorem 5.28(ii)]{nicolatilli_norm}), that are
	\begin{equation*}
		F_p(z) = e^{i\theta} \lambda\left(1-\left| \frac{z-z_0}{z-\overline{z_0}}\right|^2 \right)^{-\alpha_p}
	\end{equation*}
	where $z \in \C_+$, $\lambda = A(4\pi\sigma_p)^{-1/p}$, for some $\theta \in \R$, $z_0 \in \C_+$, may have an $L^q$ norm less or equal than $B$. Using the expression of their distribution function (see \cite[Equation (5.17)]{nicolatilli_norm}) it is immediate to obtain the following
	\begin{align*}
		\|F_p\|_q^q &= q \int_0^{\infty} t^{q-1} \nu(\{|F_p|>t\}) \, dt \\
		&= 4 \pi q \lambda^q \int_0^1 (s^{q-1-\alpha_p}-s^{q-1}) \, ds.						
	\end{align*}
	This integral is not finite for $q \leq \alpha_p$, while for $q > \alpha_p$ we have
	\begin{equation*}
		\|F_p\|_q = (4\pi)^{\frac{1}{q}-\frac{1}{p}}\left(\frac{\alpha_p}{q-\alpha_p}\right)^{\frac{1}{q}}(\sigma_p)^{-\frac{1}{p}}A,
	\end{equation*}
	which is less or equal than $B$ if and only if
	\begin{equation*}
		\frac{B}{A} \geq (4 \pi)^{\frac{1}{q} - \frac{1}{p}} \left(\dfrac{\alpha_p}{q-\alpha_p}\right)^{\frac{1}{q}} (\sigma_p)^{-\frac{1}{p}} = r_2.
	\end{equation*}
	So, in this regime of the parameters, $F_p$ solves also the problem with the $L^q$ constraint. Repeating the same argument with $q$ instead of $p$ gives the remaining part of the statement.
\end{IEEEproof}
The regime we just considered (that is $p>\alpha_q$ and $B/A \leq r_1$ or $q > \alpha_p$ and $B/A \geq r_2$) does not completely solve the problem. To obtain an estimate in the remaining regime, we have to look at the variational problem
\begin{equation}\label{eq:variational problem}
	\sup_{v \in \mathcal{C}} \int_0^{\infty} G(v(t))\,dt,
\end{equation}
where $\mathcal{C}$ is the class of all distribution functions coming from weights satisfying \eqref{eq : constraints wavelet case} and is defined as
\begin{equation}
	\begin{split}
		\mathcal{C} = \{ &v \colon (0,\infty) \to [0,\infty) \text{ nonincreasing s.t. }\\
		&p\int_0^{\infty}t^{p-1}v(t)\,dt \leq A^p,\\
		&q\int_0^{\infty}t^{q-1}v(t)\,dt \leq B^q\}.
	\end{split}
\end{equation}
We point out that this variational problem is similar to the one that has been solved in \cite{riccardi_new_optimal}. Indeed, all the information about the original problem are encoded in the function $G$ and the measure $\nu$. Our result is based on the following proposition, which is quite general and is independent of the explicit expression of the function $G$. As a matter of fact, the only assumptions on $G$ are the following:
\begin{itemize}
	\item $G$ is of class $C^1$;
	\item $G$ is strictly increasing;
	\item $G(s) \leq \min\{1,Ks\}$ for every $s>0$ and some $K>0$ (in this case $K=\beta/2\pi$).
\end{itemize}
\begin{prop}\label{prop:variational problem}
	The variational problem \eqref{eq:variational problem} has a unique maximizer, that is given by
	\begin{equation*}
		u(t) = 4\pi \max\{(\lambda_1 t^{p-1}+\lambda_2 t^{q-1})^{-\frac{1}{2\beta+1}}-1,0\},
	\end{equation*}
	for $t \in (0,\infty)$, where $\lambda_1,\lambda_2 > 0$ are uniquely determined by
	\begin{equation*}
		\begin{cases}
			p \int_0^{+\infty} t^{p-1} u(t) \,dt = A^p,\\
			q \int_0^{+\infty} t^{q-1} u(t) \,dt = B^q.
		\end{cases}
	\end{equation*}
\end{prop}
\begin{IEEEproof}[Sketch of the proof]
	We use some arguments similar to those already appeared in \cite{nicolatilli_norm, riccardi_new_optimal}. Hence, we only give a sketch of the proof.
	\begin{itemize}
		\item The first thing to prove is the existence of a maximizer. This is done using the direct method of calculus of variations (compactness is given by Helly's selection theorem).
		\item Then, one needs to enlarge $\mathcal{C}$ by removing the monotonicity assumption, while proving that maximizers are still in $\mathcal{C}$. This is done using decreasing rearrangements.
		\item Since \eqref{eq:variational problem} is a constrained optimization problem it is easy to see that if $u$ is a maximizer then
		\begin{equation*}
			G'(u(t)) = \lambda_1 t^{p-1} + \lambda_2 t^{q-1},
		\end{equation*}
		where $t \in (0,M)$, for some $M > 0$ and some Lagrange multipliers $\lambda_1,\lambda_2 \in \R$. Using the explicit expression of $G$, up to renaming the multipliers this implies
		\begin{equation*}\label{eq:expression of u}
			u(t) = \begin{cases}
				4 \pi [(\lambda_1 t^{p-1} + \lambda_2 t^{q-1})^{-\frac{1}{2\beta+1}}-1],\ &t \in (0,M)\\
				0, &t \in (M,\infty)
			\end{cases}
		\end{equation*}
		\item Using this explicit expression, one can prove that maximizers must achieve equality in both constraints, that $\lambda_1$ and  $\lambda_2$ are both positive  and that maximizers are continuous, which implies $M=T$, where $T$ is the unique positive solution of the equation $\lambda_1 T^{p-1} + \lambda_2 T^{q-1} = 1$.
		\item Finally, using the particular regime of the parameters $p$, $q$, $A$ and $B$, one proves that $\lambda_1$ and $\lambda_2$ are uniquely determined by the condition that $u$ must achieve equality in the constraints. If we interpret $u$ as a function of $t$, $\lambda_1$ and $\lambda_2$, this is done by studying the level sets of the functions
		\begin{equation*}
			\begin{split}
				f(\lambda_1, \lambda_2) &= p\int_0^T t^{p-1}u(t,\lambda_1,\lambda_2)\,dt\\
				g(\lambda_1, \lambda_2) &= q\int_0^T t^{q-1}u(t,\lambda_1,\lambda_2)\,dt.
			\end{split}
		\end{equation*}
		Indeed, one can prove that the level sets $\{f = A^p\}$ and $\{g=B^q\}$ do intersect (due to the regime of the parameters) and that the intersection is unique. This, in the end, implies the uniqueness of the maximizer.
	\end{itemize}
\end{IEEEproof}
Thanks to this proposition we can conclude the proof of Theorem \ref{th:main theorem wavelet case}.
\begin{IEEEproof}[Proof of Theorem \ref{th:main theorem wavelet case}(2)]
	Combining Theorem \ref{th:estimate norm} and Proposition \ref{prop:variational problem} we have that
	\begin{equation*}
		\|L_{F,\beta}\| \leq \int_0^{\infty} G(v(t)) \, dt \leq \int_0^T G(u(t)) \, dt,
	\end{equation*}
	which gives the optimal estimate.
	
	For the characterization of optimal weights, we notice that, due to Theorem \ref{th:estimate norm}, equality in the first estimate is achieved if and only if
	\begin{equation*}
		F(z) = e^{i\theta}\rho(|z-z_0|^2/|z-\overline{z_0}|^2),
	\end{equation*}
	for some $\theta \in \R$, $z_0 \in \C_+$ and $\rho \colon [0,1) \to [0,\infty)$ nonincreasing, while, due to Proposition \ref{prop:variational problem}, equality in the second estimate is achieved if and only if $v$, that is the distribution function of $|F|$, is equal to $u$. This allows us to reconstruct (at least, in an implicit way) the expression of optimal weights. Indeed, the super-level set $\{|F|>t\}$ is an hyperbolic disc with centre $z_0$ and $\nu$-measure equal to $u(t)$, that is the subset of $\C_+$ given by
	\begin{equation*}
		\left| \frac{z-z_0}{z-\overline{z_0}}\right|^2 < 1 - \left(1+\frac{u(t)}{4\pi}\right)^{-1}.
	\end{equation*}
	For simplicity, we let $d(z,z_0) = |z-z_0|^2/|z-\overline{z_0}|^2$.
	Clearly, for all the points on the boundary of this disk it holds
	\begin{equation*}
		t = \rho(d(z,z_0))
	\end{equation*}
	and 
	\begin{equation*}
		d(z,z_0) = 1 - \left(1+\frac{u(t)}{4\pi}\right)^{-1}.
	\end{equation*}
	Inverting the last equation for $t<T$ leads to
	\begin{equation*}
		\frac{d(z,z_0)}{1-d(z,z_0)} = u(t) = (\lambda_1 t^{p-1} + \lambda_2 t^{q-1})^{-\frac{1}{2\beta+1}}
	\end{equation*}
	and, recalling that $\psi$ denotes the inverse function of the right-hand side, we have
	\begin{equation*}
		t = \psi\left(\frac{d(z,z_0)}{1-d(z,z_0)}\right),
	\end{equation*}
	which implies
	\begin{equation*}
		\rho(d(z,z_0)) = \psi\left(\frac{d(z,z_0)}{1-d(z,z_0)}\right).
	\end{equation*}
	In conclusion, we obtained that optimal weight functions are given by
	\begin{equation*}
		F(z) = e^{i\theta} \psi\left(\frac{|z-z_0|^2/|z-\overline{z_0}|^2}{1-|z-z_0|^2/|z-\overline{z_0}|^2}\right),
	\end{equation*}
	which concludes the proof.
\end{IEEEproof}
We conclude with some remarks.
\begin{remark}
	The proof of Theorem \ref{th:main theorem wavelet case} strongly relies on the concentration estimate proved in \cite{ramos-tilli}, that is
	\begin{equation*}
		\int_{\Omega} |\W_{\psi_{\beta}}f(x,y)|^2 \, d\nu(x,y) \leq 1-\left(1+\frac{\nu(\Omega)}{4\pi}\right)^{-2\beta}
	\end{equation*}
	for every normalized $f \in H^2(\R)$ and every $\Omega \subset \C_+$ with finite $\nu$-measure, along with the characterization of optimal sets. Once this result is known in some other setting, the proof of the optimal estimate for localization operators follows the same path since, as already noticed, the variational problem \eqref{eq:variational problem} does not depend on the original setting (e.g. short-time Fourier transform, wavelet etc.). As an example, using \cite[Theorem 5.2]{ortega} one can prove analogous estimates in the ``spherical setting'' (which corresponds to spaces of holomorphic polynomials on the Riemann sphere with the Fubini-Study metric).
\end{remark}
\begin{remark}
	For one of the Lebesgue exponents, say $q$, that goes to $\infty$, one would expect to recover the result from \cite{nicolatilli_norm}, where the problem with an $L^p$ and $L^{\infty}$ constraint is addressed. For the expression of the regimes, this is true. Indeed, if $q$ is sufficiently large it is easy to see that the only condition that matters in differentiating between Theorem \ref{th:main theorem wavelet case}(1) and Theorem \ref{th:main theorem wavelet case}(2) is if $B/A$ is larger or smaller than $r_2$. Moreover, for $q \to \infty$ we have $r_2 \to (4\pi\sigma_p)^{-1/p}$, which is exactly the threshold that appears in \cite[Theorem 5.2]{nicolatilli_norm}. However, the study of the limit of the expression of the maximizer $u$ seems to be quite involved, since the quantities $\lambda_1$, $\lambda_2$ and $T$ appearing in \eqref{eq:expression of maximizer distribution function} are all defined in an implicit way. In addition to that, obtaining the estimate for the $L^p \cap L^{\infty}$ case with this limit argument would not automatically yield the uniqueness of the extremizers.
\end{remark}
\begin{remark}
	Theorem \ref{th:main theorem wavelet case} has some immediate yet interesting corollaries, which we do not state explicitly for the sake of brevity. For example using a duality argument one can immediately prove and optimal estimate for $\| | \W_{\psi_{\beta}} f|^2 \|_{L^p+L^q}$ (see \cite[Corollary 3.2]{riccardi_new_optimal} for the analogous result for the short-time Fourier transform). Moreover, thanks to the Bergman transform, one can transfer these results to functions in the Bergman space and Toeplitz operators on these spaces.
\end{remark}

\end{document}